\documentstyle[pb-diagram]{article}

\newcommand{\es}{\emptyset}

\newcommand{\ba}{\begin{array}}
\newcommand{\ea}{\end{array}}

\newtheorem{theorem}{Theorem}
\newtheorem{proposition}[theorem]{Proposition}
\newtheorem{lemma}[theorem]{Lemma}

\newcommand{\be}{\begin{enumerate}}
\newcommand{\ee}{\end{enumerate}}
\newcommand{\bi}{\begin{itemize}}
\newcommand{\ei}{\end{itemize}}
\newcommand{\bd}{\begin{description}}
\newcommand{\ed}{\end{description}}

\newcommand{\beq}{\begin{eqnarray*}}
\newcommand{\eeq}{\end{eqnarray*}}
\newcommand{\seq}{\Rightarrow}

\author{ {F.Parlamento, F.Previale }
\\Department of Mathematics,  Computer Science and Physics
\\University of Udine,  via  Delle Scienze 206, 33100 Udine, Italy.
\\Department of Mathematics
\\University of Turin, via Carlo Alberto 10, 10123 Torino, Italy
\\e-mail: {\em franco.parlamento$@$uniud.it} 
}
 \title{A simplified version of  the Sequent Calculus ${\bf G3[mic]}^=$}
\date{}

\begin{document}
\maketitle




\begin{abstract}
We  show that the Replacement Rule in the sequent calculus ${\bf G3[mic]}^= $, for first order languages with function symbols and equality, can be replaced by  the simpler
rule in which the transformed formula is not repeated in the premiss.

\end{abstract}


\


\section{Introduction}

As it is well known, the  modern  ${\bf G3}$ sequent calculi free of structural rules have evolved from the work of Gentzen through Ketonen, Kleene, Dragalin and Troelstra, until, in the words of \cite{NvP01}, ``a gem emerged".  In particular such systems  were obtained  by showing that the repetition of the principal formula in the premis(ses) of the logical rules, that Kleene, in \cite{K52}, had proposed in all reasonable cases, in most of them, could actually be dispensed with.
An extension of such calculi to logic with equality was proposed by Negri and van Plato in \cite{NvP98}. As the authors write in \cite{NvP11}, Troelstra appreciated so much their proposal, that it  was  adopted in the second edition \cite{TS00} of \cite{TS96}, thus replacing the standard axiomatic treatment of equality of the first edition.
The rules for equality proposed in \cite{NvP98} (see also \cite{NvP01}) were  the following:

\[
\ba{clcl}
 a=a, \Gamma \seq \Delta&\vbox to 0pt{\hbox{Ref}}~~~~~~~~~~~& a=b, P[v/a], P[v/b],  \Gamma \seq \Delta &\vbox to 0pt{\hbox{Repl}}\\
\cline{1-1}\cline{3-3}
\Gamma\seq \Delta&&a=b, P[v/a], \Gamma  \seq \Delta
\ea
\]
where $\Gamma$ and $\Delta$ are finite multisets of formulae, $P$ is an atomic formula and
$P[v/a]$ and $P[v/b]$ denote  the result of the substitution of the variable $v$ by the individual constants $a$ and $b$ respectively.

The rules   adopted in \cite{TS00} are the above $\hbox{Ref}$ and $\hbox{Repl}$, with $a$ and $b$ replaced by arbitrary terms $s$ and $r$ and  the proviso that   $v$ does not occur in $s$ and $r$.

Our purpose is to show that, in line with the evolution from Kleene's rules to the present ones, also the repetition of the formula $P[v/s]$ in the premiss of the rule $\hbox{Repl}$, meant  to ensure the admissibility of the contraction rule,  can actually be dispensed with.

More precisely, we will refer  to the  multisuccedent systems ${\bf G3[mic]}^=$, for  minimal, intuitionistic and classical logic,   in \cite{N16}, with  the rules $\hbox{Ref}$  and $\hbox{Repl}$  extended to arbitrary terms, and will show that
 in ${\bf G3[mic]}^=$ the rule $\hbox{Repl}$ can be replaced by the, apparently weaker, rule $\hbox{Repl}^-$:
\[
\ba{c}
s=r, P[v/r], \Gamma\seq \Delta\\
\cline{1-1}
 s=r, P[v/s], \Gamma\seq \Delta
\ea
\]
by showing  that  the  rule $\hbox{Repl}$ is admissible in  the  systems ${{\bf G3[mic]}^{=}}^-$  obtained by replacing in ${\bf G3[mic]}^{=}$ the rule $\hbox{Repl}$ by $\hbox{Repl}^-$. As a consequence ${\bf G3[mic]}^=$   and ${{\bf G3[mic]}^=}^-$  are equivalent.

\ 

\section{Admissibility of $\hbox{Repl}$ in ${{\bf G3[mic]}^=}^-$ }
Let $\hbox{Repl}_1^-$ and $\hbox{Repl}_1$ be the rules $\hbox{Repl}^-$ and $\hbox{Repl}$ 
as represented in the Introduction, in which it is required that there is exactly one occurrence of $v$ in $A$, i.e. only one occurrence of $r$ is replaced by $s$.

\begin{lemma}\label{single}
 $\hbox{Repl}^-$  is derivable from  $\hbox{Repl}_1^-$.
\end{lemma}

{\bf Proof} By induction on the number of occurrences of $v$ in $A$.
If such a number  is $n+1$, with $n\geq 1$, then let $A'$ be obtained by replacing one occurrence of $v$ in $A$ by a new variable $v'$.
$A[v/r]$ coincides with $(A'[v/r])[v'/r]$. Thus from $$ s=r, A[v/r], \Gamma\seq \Delta$$ by  $\hbox{Repl}_1^-$ 
we can obtain $$ s=r, (A'[v/r])[v'/s], \Gamma\seq \Delta$$ that coincides with
 $$s=r, (A'[v'/s])[v/r], \Gamma \seq \Delta$$ Since  there are $n$ occurrences of $v$ in $A'[v'/s]$,
 by the induction hypothesis, from the latter sequent, using $\hbox{Repl}_1^-$ we can derive
 $$ s=r, (A'[v'/s])[v/s], \Gamma \seq \Delta$$ that coincides with $$ s=r, A[v/s], \Gamma \seq \Delta$$ $\Box$

\begin{lemma}\label{single2}
 $\hbox{Repl}$ is derivable  from  $\hbox{Repl}_1$ and  the left weakening rule.
\end{lemma}

{\bf Proof} As in the proof of Lemma \ref{single}  let $A'$ be obtained by replacing one of the $n+1$ occurrences of $v$ in $A$ by a new variable $v'$.
The premiss
 $$s=r, A[v/s], A[v/r], \Gamma  \seq \Delta$$ of $\hbox{Repl}$
   coincides with
  $$ s=r, (A'[v/s])[v'/s], (A'[v/r])[v'/r], \Gamma \seq \Delta$$
from which by the left weakening rule 
 we obtain 
$$s=r, (A'[v/s])[v'/s], (A'[v/r])[v'/s], (A'[v/r])[v'/r], \Gamma \seq \Delta$$

   Then an application of $\hbox{Repl}_1$ yields
   $$ s=r, (A'[v/s])[v'/s], (A'[v/r])[v'/s], \Gamma \seq \Delta$$
   namely
    $$s=r, (A'[v'/s])[v/s], (A'[v'/s])[v/r], \Gamma \seq \Delta$$
    from which, by the induction hypothesis,
    we can derive 
    $$ s=r, (A'[v'/s])[v/s], \Gamma \seq \Delta$$ i.e. 
    $$ s=r, A[v/s], \Gamma \seq \Delta$$
   $\Box$

\

Let ${{\bf G3[mic]}_1^=}^-$ be  obtained by replacing $\hbox{Repl}^-$ by $\hbox{Repl}_1^-$ in ${{\bf G3[mic]}^=}^-$ .

For the sake of notational brevity in the following we will denote ${{\bf G3[mic]}^=}^-$  and  ${{\bf G3[mic]}_1^=}^-$ also by $S$ and $S_1$ respectively.

\begin{lemma}\label{weakening}
The weakening rules are height preserving admissible in  $S$  and $S_1$, i.e.  if $\Gamma\seq \Delta$ has a derivation in $S$ ($S_1$)  of  height $\leq h$, 

then also $F, \Gamma \seq \Delta$ and $\Gamma\seq \Delta, F$ have a derivation in $S$ ($S_1$) of height $\leq h$.
 \end{lemma}
 
 \begin{lemma}
 \bi
 
\item[a)] Derivability in $S_1$ of $Contr^=$ 

 $s=r, \Gamma\seq \Delta$ is derivable in  $S_1$ from  $ s=r, s=r, \Gamma \seq \Delta$

\item[b)]   Admissibility  in  $S_1$ of $Symm$:

 In  $S_1$ and the left weakening rule,   $ r=s, \Gamma \seq \Delta$  is derivable from
 
   $s=r, \Gamma \seq \Delta$.  The same holds for  $S$.
\ei

 \end{lemma}

{\bf Proof}
$a)$  Since $s=r$ coincides with $A[v/r]$, where $A$ is $s=v$, the following is a derivation in $S$ of $s=r, \Gamma \seq \Delta$ from $s=r, s=r,\Gamma \seq \Delta$

\[
\ba{cl}
{\cal D}&\\
s=r, s=r, \Gamma \seq \Delta&\vbox to 0pt{$\hbox{Repl}_1^-$}\\
\cline{1-1}
s=r, s=s, \Gamma \seq \Delta &\vbox to 0pt{$\hbox{Ref}$}\\
\cline{1-1}
s=r, \Gamma \seq \Delta &
\ea
\]

$b)$ The following is a  derivation in $S_1$ and the left weakening rule of 

$r=s, \Gamma \seq \Delta$ from  $s=r, \Gamma\seq \Delta$:
\[
\ba{cl}
 s=r,  \Gamma\seq \Delta&\\
 \cline{1-1}
 s=r, r=r, \Gamma\seq \Delta&\vbox to 0pt{$\hbox{Repl}_1^-$}\\
\cline{1-1}
 s=r, r=s, \Gamma \seq \Delta&\vbox to 0pt{$\hbox{Repl}_1^-$}\\
\cline{1-1}
r=r, r=s, \Gamma \seq \Delta&\vbox to 0pt{$\hbox{Ref}$}\\
\cline{1-1}
 r=s, \Gamma \seq \Delta
\ea
\]

$\Box$

\

\begin{proposition}\label{admissibility}
$\hbox{Repl}$ is admissible in   ${{\bf G3[mic]}^=}^-$
\end{proposition}

{\bf Proof~} We have to show that the applications of the rule $\hbox{Repl}$ can be eliminated from the derivations in  $S +\hbox{Repl}$. Since the left weakening rule is admissible in $S$, by Lemma \ref{single2},
we can transform a derivation in  $S  + \hbox{Repl}$ into a derivation with  the same endsequent in  $S +\hbox{Repl}_1$. Thus by Lemma \ref{single}   it suffices  to show that a derivation ${\cal D}$ in $S_1$ of 

$$s=r, A[v/s], A[v/r], \Gamma \seq \Delta$$ with $v$ that does not occur in $s,r$ and has a single occurrence in $A$, can be transformed into a derivation ${\cal D}'$  in $S_1$ of $ s=r, A[v/s], \Gamma \seq \Delta$.
The proof is by induction on the height of derivations, but for the induction argument to go through we have to generalize the statement to be proved.
In fact,  assume that $A[v/s]$ has the form $A^\circ[u/q, v/s]$ and the given derivation of $s=r, A[v/s], A[v/r], \Gamma \seq \Delta$ has the form:

\[
\ba{c}
{\cal D}_0\\
q=p, s=r, A^\circ[u/q, v/s], A^\circ[u/p,v/r], \Gamma\seq \Delta\\
\cline{1-1}
q=p, s=r, A^\circ[u/q, v/s], A^\circ[u/q,v/r], \Gamma\seq \Delta
\ea
\]
Then since $A^\circ[u/q, v/s]$ and $ A^\circ[u/p,v/r]$ do not have the form $B[v/s]$ and $B[v/r]$ we could not  apply the induction hypothesis to 
${\cal D}_0$.

To overcame that problem,  we generalize the statement to be proved as follows.
Let $\vec{q}$ and  $\vec{p}$ be the sequences of terms 
  $q_1,\ldots q_n$  and $p_1,\ldots p_n$  and, similarly, let  $\vec{u}$ stand for the sequence of variables $u_1,\ldots u_n$ assumed to be distinct from one another and from $v$ and not occurring in $\vec{q},\vec{p}, s,r$.   
$\vec{q}=\vec{p}$ stands for the sequence of equalities $q_1=p_1, \ldots, q_n=p_n$ and $[\vec{u}/\vec{q}]$ for the substitution $[u_1/q_1, \ldots, u_n/q_n]$ and similarly for $[\vec{u}/\vec{p}]$.  We proceed by induction on the height $h({\cal D})$ of ${\cal D}$ to show that if   ${\cal D}$ is a derivation in $S_1$ of

$$
\vec{q}=\vec{p}, s=r, A[\vec{u}/\vec{q},  v/s],  A[\vec{u}/\vec{p},  v/r ], \Gamma \seq \Delta
$$
where each one of the variables in $\vec{u}$ and $v$ has exactly one occurrence in $A$, 
then ${\cal D}$ can be transformed into a derivation ${\cal D}'$ in $S_1$ of 
$$
 \vec{q}=\vec{p}, s=r, A[\vec{u}/\vec{q},  v/s], \Gamma \seq \Delta
$$

The statement we are actually interested in, that yields the admissibility of $\hbox{Repl}_1$, therefore of   $\hbox{Repl}$, in $S$,  follows by letting $n=0$.

If $h({\cal D})=0$, then ${\cal D}$ reduces to a logical axiom. 
If $\Gamma\cap \Delta \neq \es$ or one of $ \vec{q}=\vec{p}, s=r, A[\vec{u}/\vec{q},  v/s]$  belongs to $\Delta$,  then 
also $\vec{q}=\vec{p}, s=r, A[\vec{u}/\vec{q},  v/s] ], \Gamma \seq \Delta$ is an axiom and we are done.
Otherwise $A[\vec{u}/\vec{p}, v/r]\in \Delta$. But  
then also $ \vec{q}=\vec{p}, s=r, A[\vec{u}/\vec{p}, v/r], \Gamma \seq \Delta$ is an axiom and as ${\cal D}'$ we can take:

\[
\ba{cl}
\vec{q}=\vec{p}, s=r, A[\vec{u}/\vec{p}, v/r], \Gamma \seq \Delta&\vbox to 0pt{$\hbox{Repl}_1^-~n+1-times$ }\\
\vdots&\\
\vec{q}=\vec{p}, s=r, A[\vec{u}/\vec{q}, v/s], \Gamma  \seq \Delta&
\ea
\]
If $h({\cal D})>0$ and ${\cal D}$ ends with a logical inference (see \cite{N16} for the complete list), since $A$ is atomic, neither $A[\vec{u}/\vec{q},v/s]$ nor $A[\vec{u}/\vec{p}, v/r]$ nor any of $\vec{q}=\vec{p}$ and $s=r$ can be the principal formula of such an inference and the conclusion follows immediately from the induction hypothesis. The same applies if ${\cal D}$ ends with a $\hbox{Ref}$-inference.
If ${\cal D}$ ends with a $\hbox{Repl}_1^-$-inference we distinguish the following cases.

Case 1. The last inference of ${\cal D}$ does not introduce any of the shown occurrences of $\vec{q}$, $\vec{p}$, $s$ or $r$.

Case 1.1 $A$ is of the form $A^\circ[u/q]$ and  ${\cal D}$ has  the form:

\[
\ba{c}
{\cal D}_0\\
 q=p, \vec{q}=\vec{p}, s=r, A^\circ[u/q,\vec{u}/\vec{q},  v/s], A^\circ[u/p,\vec{u}/\vec{p}, v/r], \Gamma' \seq \Delta\\
\cline{1-1}
 q=p,  \vec{q}=\vec{p}, s=r, A^\circ [u/q, \vec{u}/\vec{q}, v/s],  A^\circ[u/q,\vec{u}/\vec{p}, v/r], \Gamma' \seq \Delta
\ea
\]

By the induction hypothesis applied to ${\cal D}_0$ and  $\Gamma'$,  $q=p, \vec{q}=\vec{p}$ and $A^\circ $ in place of $\Gamma$ , $\vec{q}=\vec{p}$ and $A$ respectively, we have a derivation ${\cal D}'_0$ in $S_1$ of 
$$q=p, \vec{q}=\vec{p}, s=r, A^\circ[u/q,\vec{u}/\vec{q},  v/s], \Gamma'  \seq \Delta$$
that can be taken as  ${\cal D}'$.

Case 1.2  $A$ is of the form $A^\circ[u/q]$ and  ${\cal D}$ has  the form:

\[
\ba{c}
{\cal D}_0\\
 q=p, \vec{q}=\vec{p}, s=r, A^\circ[u/p,\vec{u}/\vec{q},  v/s], A^\circ[u/q,\vec{u}/\vec{p}, v/r], \Gamma' \seq \Delta\\
\cline{1-1}
 q=p,  \vec{q}=\vec{p}, s=r, A^\circ [u/q, \vec{u}/\vec{q}, v/s],  A^\circ[u/q,\vec{u}/\vec{p}, v/r], \Gamma' \seq \Delta
\ea
\]
By height-preserving weakening we have a derivation ${\cal D}_0^w$ of the same height as ${\cal D}_0$ of
\[
p=q,  q=p,  \vec{q}=\vec{p}, s=r, A^\circ [u/p, \vec{u}/\vec{q}, v/s],  A^\circ[u/q,\vec{u}/\vec{p}, v/r], \Gamma' \seq \Delta
 \]
By induction hypothesis there is a derivation  ${\cal D}_0^{w'}$ in $S_1$ of 

\[
p=q,  q=p,  \vec{q}=\vec{p}, s=r, A^\circ [u/p, \vec{u}/\vec{q}, v/s], \Gamma' \seq \Delta
 \]
Then ${\cal D}'$  can be obtained from   the following derivation in $S_1+\hbox{Contr}^= + \hbox{Symm}$,
 thanks to the derivability in $S_1$ of $\hbox{Contr}^=$ and the admissibility in $S_1$ of $\hbox{Symm}$:

\[
\ba{cl}
p=q,  q=p,  \vec{q}=\vec{p}, s=r, A^\circ [u/p, \vec{u}/\vec{q}, v/s], \Gamma' \seq \Delta&\vbox to 0pt{$\hbox{Repl}_1^-$ }\\
\cline{1-1}
p=q,  q=p,  \vec{q}=\vec{p}, s=r, A^\circ [u/q, \vec{u}/\vec{q}, v/s], \Gamma' \seq \Delta&\vbox to 0pt{$\hbox{Symm}$ }\\
\cline{1-1}
 q=p, q=p,  \vec{q}=\vec{p}, s=r, A^\circ [u/q, \vec{u}/\vec{q}, v/s], \Gamma' \seq \Delta&\vbox to 0pt{$\hbox{Contr}^=$ }\\
 \cline{1-1}
 q=p,  \vec{q}=\vec{p}, s=r, A^\circ [u/q, \vec{u}/\vec{q}, v/s], \Gamma' \seq \Delta&
 \ea
\]

Case 1.3   The last $\hbox{Repl}_1^-$ inference of ${\cal D}$ acts by means of $q=p$ inside $\Gamma'$. Then the conclusion follows by applying the induction hypothesis and then the same $\hbox{Repl}_1^-$-inference .

\

Case 2.  The last inference of ${\cal D}$ does  introduce one of the  shown occurrences of $\vec{q}$, $\vec{p}$, $s$ or $r$. Without loss of generality we may assume that is either  $r$ or $s$.

\

Case 2.1  The  last inference of ${\cal D}$ introduces  $r$.

Case 2.1.1 $A$ has the form $A^\circ[u/q]$, $v$ occurs in $q$ and ${\cal D}$ has  the form:

\[
\ba{c}
{\cal D}_0\\
q[v/r]=p,~ \vec{q}=\vec{p},  ~s=r,  A^\circ[\vec{u}/\vec{q},  ~ u/q[v/s]],  A^\circ[\vec{u}/\vec{p},~ u/p ], \Gamma' \seq \Delta\\
\cline{1-1}
 q[v/r]=p, ~ \vec{q}=\vec{p}, ~s=r,  A^\circ[\vec{u}/\vec{q}, ~ u/q[v/s]],  A^\circ[\vec{u}/\vec{p}, ~ u/q[v/r]], \Gamma'  \seq \Delta
\ea
\]
By height-preserving weakening we have a derivation ${\cal D}_0^w$ of the same height as ${\cal D}_0$ of

\[
 q[v/r]=p,~ q[v/s]=p,~\vec{q}=\vec{p},~s=r,  A^\circ[\vec{u}/\vec{q}, ~ u/q[v/s]],  A^\circ [\vec{u}/\vec{p}, ~ u/p ], \Gamma' \seq \Delta
\]
By induction hypothesis there is a derivation ${\cal D}_0^{w'}$ in $S_1$ of 

\[
q[v/r]=p,~ q[v/s]=p, \vec{q}=\vec{p},~s=r,  A^\circ[\vec{u}/\vec{q},  ~ u/q[v/s]], \Gamma' \seq \Delta
\]
 Then ${\cal D}'$  can be obtained from:
\[
\ba{cl}
{\cal D}_0^{w'}&\\
q[v/r]=p,~ q[v/s]=p, ~\vec{q}=\vec{p},~s=r,  A^\circ[\vec{u}/\vec{q}, ~ u/q[v/s]], \Gamma' \seq \Delta
&\vbox to 0pt{$\hbox{Symm}$ }\\
\cline{1-1}
 q[v/r]=p,~ q[v/s]=p, ~\vec{q}=\vec{p}, ~r=s,  A^\circ[\vec{u}/\vec{q}, ~ u/q[v/s]], \Gamma'\seq \Delta&\vbox to 0pt{$\hbox{Repl}_1^-$ }\\
\cline{1-1}
q[v/r]=p,~ q[v/r]=p, ~\vec{q}=\vec{p}, ~r=s,  A^\circ[\vec{u}/\vec{q}, ~ u/q[v/s]], \Gamma' \seq \Delta&\vbox to 0pt{$\hbox{Contr}^=$ }\\
\cline{1-1}
q[v/r]=p,~ \vec{q}=\vec{p},~r=s,  A^\circ[\vec{u}/\vec{q},~ u/q[v/s]], \Gamma' \seq \Delta&\vbox to 0pt{$\hbox{Symm}$ }\\
\overline{ q[v/r]=p, ~ \vec{q}=\vec{p},~s=r,  A^\circ[\vec{u}/\vec{q}, ~ u/q[v/s]], \Gamma'\seq \Delta}
\ea
\]

Case 2.1.2  $r$ is of the form $r^\circ[u/q]$ and ${\cal D}$ has the form:

\[
\ba{c}
{\cal D}_0\\
q=p,~\vec{q}=\vec{p},  ~ s=r^\circ[u/q],  A[\vec{u}/\vec{q}, ~ v/s], A[\vec{u}/\vec{p},  v/r^\circ[u/p]], \Gamma' \seq \Delta\\
 \cline{1-1}
 q=p, ~\vec{q}=\vec{p},  ~ s=r^\circ[u/q],  A[\vec{u}/\vec{q}, ~ v/s], A[\vec{u}/\vec{p}, ~ v/r^\circ[u/q]],  \Gamma' \seq \Delta
\ea
\]
By height-preserving weakening we have a derivation ${\cal D}_0^w$ of the same height as ${\cal D}_0$ of

\[
q=p,~ \vec{q}=\vec{p},  ~s=r^\circ[u/q], ~s=r^\circ[u/p],   A[\vec{u}/\vec{q},~v/s], A[\vec{u}/\vec{p}, ~ v/r^\circ[u/p]], \Gamma' \seq \Delta
\]

By induction hypothesis there is a derivation ${\cal D}_0^{w'}$  in $S_1$ of 
\[
q=p,~ \vec{q}=\vec{p}, ~ s=r^\circ[u/q],~s=r^\circ[u/p],   A[\vec{u}/\vec{q}, ~ v/s], \Gamma' \seq \Delta
\]
Then  ${\cal D}'$  can  be obtained from:
\[
\ba{cl}
{\cal D}_0^{w'}&\\
q=p, ~\vec{q}=\vec{p},  ~ s=r^\circ[u/q],~s=r^\circ[u/p],   A[\vec{u}/\vec{q}, ~ v/s], \Gamma \seq \Delta&\vbox to 0pt{$\hbox{Repl}_1^-$ }\\
\cline{1-1}
q=p,~\vec{q}=\vec{p},  ~ s=r^\circ[u/q],~s=r^\circ[u/q],   A[\vec{u}/\vec{q}, ~ v/s], \Gamma \seq \Delta&\vbox to 0pt{$\hbox{Contr}^=$ }\\
\overline{       q=p, ~\vec{q}=\vec{p}, ~s=r^\circ[u/q],   A[\vec{u}/\vec{q}, ~ v/s],  \Gamma\seq \Delta      }&
\ea
\]

Case 2.1.3  $r$ is of the form $r^\circ[u/q]$ and ${\cal D}$ has the form:

\[
\ba{c}
{\cal D}_0\\
q=p,~\vec{q}=\vec{p},  ~ s=r^\circ[u/p],  A[\vec{u}/\vec{q}, ~ v/s], A[\vec{u}/\vec{p}, ~ v/r^\circ[u/q]]], \Gamma'  \seq \Delta\\
 \cline{1-1}
q=p, ~\vec{q}=\vec{p},  ~ s=r^\circ[u/q],  A[\vec{u}/\vec{q}, ~ v/s], A[\vec{u}/\vec{p}, ~ v/r^\circ[u/q]],  \Gamma'  \seq \Delta
\ea
\]

By height-preserving weakening we have a derivation ${\cal D}_0^w$ in $S_1$ of the same height as ${\cal D}_0$ of

\[
 q=p,~ \vec{q}=\vec{p},  ~s=r^\circ[u/p], ~s=r^\circ[u/q],   A[\vec{u}/\vec{q}, ~ v/s], A[\vec{u}/\vec{p}, ~ v/r^\circ[u/q]], \Gamma' \seq \Delta
\]
By induction hypothesis there is a derivation ${\cal D}_0^{w'}$ in $S_1$ of 
\[
q=p,~ \vec{q}=\vec{p}, ~ s=r^\circ[u/p],~s=r^\circ[u/q],   A[\vec{u}/\vec{q}, ~ v/s], \Gamma' \seq \Delta
\]
 Then  ${\cal D}'$  can  be obtained from:

\[
\ba{cl}
{\cal D}_0^{w'}\\
 q=p,~ \vec{q}=\vec{p}, ~ s=r^\circ[u/p],~s=r^\circ[u/q],   A[\vec{u}/\vec{q}, ~ v/s], \Gamma' \seq \Delta&\vbox to 0pt{$\hbox{Repl}_1^-$ }\\
\cline{1-1}
 q=p,~ \vec{q}=\vec{p}, ~ s=r^\circ[u/q],~s=r^\circ[u/q],   A[\vec{u}/\vec{q}, ~ v/s], \Gamma' \seq \Delta&\vbox to 0pt{$\hbox{Contr}^=$ }\\
\cline{1-1}
q=p,~ \vec{q}=\vec{p}, ~ s=r^\circ[u/q],  A[\vec{u}/\vec{q}, ~ v/s], \Gamma' \seq \Delta
\ea
\]

Case 2.2    The  last inference of ${\cal D}$ introduces  $s$.

Case 2.2.1 $A$ has the form $A^\circ[u/q]$, $v$ occurs in $q$ and ${\cal D}$ has  the form:

\[
\ba{c}
{\cal D}_0\\
 q[v/s]=p,~ \vec{q}=\vec{p},  ~s=r,  A^\circ[\vec{u}/\vec{q}, ~ u/p],  A^\circ[\vec{u}/\vec{p},~ u/q[v/r] ], \Gamma' \seq \Delta\\
\cline{1-1}
 q[v/s]=p, ~ \vec{q}=\vec{p}, ~s=r,  A^\circ[\vec{u}/\vec{q}, ~ u/q[v/s]], A^\circ[\vec{u}/\vec{p},~ u/q[v/r]], \Gamma' \seq \Delta
\ea
\]

By height-preserving weakening we have a derivation ${\cal D}_0^w$ of the same height as ${\cal D}_0$ of

\[
q[v/s]=p,~ p=q[v/r],~\vec{q}=\vec{p},~s=r,  A^\circ[\vec{u}/\vec{q}, ~ u/p],  A^\circ [\vec{u}/\vec{p},~ u/q[v/r] ], \Gamma'  \seq \Delta
\]
By induction hypothesis there is a derivation ${\cal D}_0^{w'}$ in $S_1$ of 

\[
q[v/s]=p,~ p=q[v/r], \vec{q}=\vec{p},~s=r,  A^\circ[\vec{u}/\vec{q}, ~ u/p], \Gamma'  \seq \Delta
\]
Then  ${\cal D}'$  can  be obtained from:

\[
\ba{cl}
{\cal D}_0^{w'}&\\
q[v/s]=p,~ p=q[v/r], \vec{q}=\vec{p},~s=r,  A^\circ[\vec{u}/\vec{q}, ~ u/p], \Gamma'  \seq \Delta&\vbox to 0pt{$\hbox{Repl}_1^-$ }\\
\cline{1-1}
q[v/s]=p,~ p=q[v/r], \vec{q}=\vec{p},~s=r,  A^\circ[\vec{u}/\vec{q}, ~ u/q[v/s]], \Gamma' \seq \Delta&\vbox to 0pt{$\hbox{Repl}_1^-$ }\\
\cline{1-1}
q[v/s]=p,~ p=q[v/s], \vec{q}=\vec{p},~s=r,  A^\circ[\vec{u}/\vec{q}, ~ u/q[v/s]], \Gamma'  \seq \Delta&\vbox to 0pt{$\hbox{Symm}$ }\\
\cline{1-1}
\cline{1-1}
q[v/s]=p,~ q[v/s]=p, \vec{q}=\vec{p},~r=s,  A^\circ[\vec{u}/\vec{q}, ~ u/q[v/s]], \Gamma' \seq \Delta&\vbox to 0pt{$\hbox{Contr}^=$ }\\
\cline{1-1}
q[v/s]=p, \vec{q}=\vec{p},~r=s,  A^\circ[\vec{u}/\vec{q},~ u/q[v/s]], \Gamma' \seq \Delta&
\ea
\]

Case 2.2.2  $s$ is of the form $s^\circ[u/q]$ and ${\cal D}$ has the form:

\[
\ba{c}
{\cal D}_0\\
q=p,~\vec{q}=\vec{p},  ~ s^\circ[u/q]=r,  A[\vec{u}/\vec{q},~ v/s^\circ[u/p]], A[\vec{u}/\vec{p}, ~v/r], \Gamma' \seq \Delta\\
 \cline{1-1}
q=p, ~\vec{q}=\vec{p},  ~ s^\circ[u/q]=r,  A[\vec{u}/\vec{q},~v/s^\circ[u/q], A[\vec{u}/\vec{p},~v/r],  \Gamma'  \seq \Delta
\ea
\]
By height-preserving weakening we have a derivation ${\cal D}_0^w$ of the same height as ${\cal D}_0$ of

\[
 q=p,~ \vec{q}=\vec{p},  ~s^\circ[u/q]=r, ~s^\circ[u/p]=r,   A[\vec{u}/\vec{q},~ v/s^\circ[u/p]], A[\vec{u}/\vec{p},  ~ v/r], \Gamma' \seq \Delta
\]
By induction hypothesis there is a derivation ${\cal D}_0^{w'}$ in $S_1$ of 
\[
q=p,~ \vec{q}=\vec{p}, ~  ~s^\circ[u/q]=r, ~s^\circ[u/p]=r,  A[\vec{u}/\vec{q}, ~v/s^\circ[u/p]], \Gamma' \seq \Delta
\]
 Then ${\cal D}'$ can be obtained from:

\[
\ba{cl}
{\cal D}_0^{w'}&\\
q=p,~ \vec{q}=\vec{p},  ~s^\circ[u/q]=r, ~s^\circ[u/p]=r,  A[\vec{u}/\vec{q},~v/s^\circ[u/p] ], \Gamma' \seq \Delta
&\vbox to 0pt{$\hbox{Repl}_1^-$ }\\
\cline{1-1}
q=p,~ \vec{q}=\vec{p},  ~s^\circ[u/q]=r, ~s^\circ[u/p]=r,  A[\vec{u}/\vec{q}, ~ v/s^\circ[u/q]], \Gamma' \seq \Delta
&\vbox to 0pt{$\hbox{Repl}_1^-$ }\\
\cline{1-1}
q=p,~ \vec{q}=\vec{p},  ~s^\circ[u/q]=r, ~s^\circ[u/q]=r,  A[\vec{u}/\vec{q}, ~ v/s^\circ[u/q]], \Gamma' \seq \Delta&\vbox to 0pt{$\hbox{Contr}^=$ }\\
\cline{1-1}
\overline{      q=p, ~\vec{q}=\vec{p},~s^\circ[u/q]=r,   A[\vec{u}/\vec{q}, ~ v/s^\circ[u/q]],  \Gamma \seq \Delta      }&
\ea
\]

Case 2.2.3  $s$ is of the form $s^\circ[u/q]$ and ${\cal D}$ has the form:

\[
\ba{c}
{\cal D}_0\\
q=p,~\vec{q}=\vec{p},  ~ s^\circ[u/p]=r,  A[\vec{u}/\vec{q}, ~ v/ s^\circ[u/q]], A[\vec{u}/\vec{p}, ~ v/r], \Gamma'  \seq \Delta\\
 \cline{1-1}
q=p, ~\vec{q}=\vec{p},  ~ s^\circ[u/q]= r,  A[\vec{u}/\vec{q}, ~ v/s^\circ[u/q]], A[\vec{u}/\vec{p}, ~ v/r],  \Gamma' \seq \Delta
\ea
\]

By height-preserving weakening we have a derivation ${\cal D}_0^w$ of the same height as ${\cal D}_0$ of

\[
q=p,~ \vec{q}=\vec{p},  ~s^\circ[u/p]=r, ~s^\circ[u/q]=r,   A[\vec{u}/\vec{q}, ~ v/s^\circ[u/q]], A[\vec{u}/\vec{p}, ~ v/r], \Gamma'  \seq \Delta
\]
By induction hypothesis there is a derivation ${\cal D}_0^{w'}$ in $S_1$ of 
\[
q=p,~ \vec{q}=\vec{p}, ~s^\circ[u/p]=r, ~s^\circ[u/q]=r,    A[\vec{u}/\vec{q}, ~ v/s^\circ[u/q]], \Gamma' \seq \Delta
\]
Then  ${\cal D}'$  can  be obtained from:

\[
\ba{cl}
{\cal D}_0^{w'}\\
 q=p,~ \vec{q}=\vec{p}, ~s^\circ[u/p]=r, ~s^\circ[u/q]=r,    A[\vec{u}/\vec{q}, ~ v/s^\circ[u/q]], \Gamma' \seq \Delta&\vbox to 0pt{$\hbox{Repl}_1^-$ }\\
\cline{1-1}
 q=p,~ \vec{q}=\vec{p},~s^\circ[u/q]=r,  ~s^\circ[u/q]=r,  A[\vec{u}/\vec{q}, ~ v/s^\circ[u/q]], \Gamma' \seq \Delta&\vbox to 0pt{$\hbox{Contr}^=$ }\\
\cline{1-1}
 q=p,~ \vec{q}=\vec{p},~s^\circ[u/q]=r,    A[\vec{u}/\vec{q}, ~ v/s^\circ[u/q]], \Gamma' \seq \Delta
\ea
\]$\Box$

\subsection{Equivalence between  ${\bf G3[mic]}^=$ and  ${{\bf G3[mic]}^=}^-$}

\begin{theorem}
A sequent is derivable in  ${\bf G3[mic]}^=$ if and only if it is derivable in    ${{\bf G3[mic]}^=}^-$.
\end{theorem}

{\bf Proof} Let ${\cal D}$  be a derivation in  ${\bf G3[mic]}^=$  of $\Gamma\seq \Delta$. By induction on the height of ${\cal D}$, thanks to Proposition 
\ref{admissibility},  it is straightforward that ${\cal D}$ can be transformed into a derivation in   ${{\bf G3[mic]}^=}^-$.
Conversely, given a derivation ${\cal D}$ in  ${{\bf G3[mic]}^=}^-$ of $\Gamma\seq \Delta$,
it suffices to apply the admissibility of the left weakening rule, to show that ${\cal D}$ can be transformed into a derivation in  ${\bf G3[mic]}^=$
of $\Gamma\seq \Delta$. $\Box$

\end{document}